\documentclass{article}

\usepackage{ifxetex,ifluatex}
\if\ifxetex T\else\ifluatex T\else F\fi\fi T%
  \usepackage{fontspec}
\else
  \usepackage[T1]{fontenc}
  \usepackage[utf8]{inputenc}
  \usepackage{lmodern}
\fi

\usepackage{hyperref}

\title{Interpolation by the Exact Inversion of the Gram Matrix}
\author{John J Spitzer}

\usepackage{sagetex}

\begin{document}
\maketitle

\begin{abstract}
Using a lemma of Davis on Gram matrices applied to the classical Orthogonal Polynomials to generate reproducing kernel interpolation over the classical domains for polynomials. These kernels have terms which are exact over the rational ring. 
The Condition Numbers are readily shown to get very large with the size of the Gram matrices as expected.
The calculation of the error variances for trigonometric functions and the exponential show a significant improvement over the equivalent Taylor expansion variances. 
\end{abstract}
\section{Introduction}
It is well known that Gram matrices which are Hilbert matrices or in more general Hankel matrices they can have very large condition numbers. This makes solving the normal equations very difficult in finite precision arithmetic. The root problem here is the difference of large integers.
However the classical orthogonal polynomials make diagonal Gram matrices. This can be used together with the lemma of Davis to solve the normal equations. This relies on the polynomial expansion of the classical orthogonal polynomials.
Note that the inner products for the normal equations are integrals over three different domains with the appropriate weights.
\begin{center}
\begin{tabular}{lrc}
Name & Range & Weight\\
\hline
Legendre & (1,1) & 1 \\
Laguerre & $(0,\infty)$ & $\exp{-x}$ \\
Hermite & $(-\infty,\infty)$ & $\exp{-x^2}$
\end{tabular}
\end{center}
These three domains are the usual common domains. Any finite domain is accessible through an affine transformation of the Legendre case.
Notice that the Legendre and Hermite domains are symmetric about the origin. Therefore the even and odd cases must be treated separately.
\section{The Lemma of Davis}
\newtheorem{thm}{Definition}
\begin{thm}[Davis]\label{def1}
From Davis~\cite[8.7.1]{Dav}

Given the sequence $(x_1,x_2,x_3,\ldots,x_n)$ 
of elements in an inner product space.  
We have the $n \times n$ matrix
\[G((x_i,x_j)) = \left[
\begin{array}{ccc}
(x_1,x_1) & (x_1,x_2) & \ldots,(x_1,x_n) \\ 
. & . &\ldots \\
. & . &\ldots \\
(x_n,x_1) & (x_n,x_2) & \ldots,(x_n,x_n)
\end{array}
\right]
\]
which is known as the Gram matrix of $(x_1,x_2,x_3,\ldots,x_n)$ . Its determinant
\[ \mathrm{g}(x_1,x_2,x_3,\ldots,x_n) = |(x_i,x_j)| = |(x_j,x_i)|
\]
is known as the Gram determinant of the elements.
\end{thm}

\newtheorem{lemma}{Lemma}
\begin{lemma}[Davis]\label{lm1}
According to (Davis~\cite[8.7.1]{Dav}) we have

Let
\[y_i = \sum_{j=1}^n (a_{ij}x_j)\]
Let A designate the matrix $(a_{ij})$ and $\tilde{A}$  be its conjugate transpose $(\tilde {a_{ji}})$. Then
\[\mathrm{G}(y_1,y_2,y_3,\ldots,y_n) = A\mathrm{G}(x_1,x_2,x_3,\ldots,x_n)\tilde{A}\]
and
\[\mathrm{g}(y_1,y_2,y_3,\ldots,y_n) = |det(A)|^2\mathrm{g}(x_1,x_2,x_3,\ldots,x_n)\]

\end{lemma}
\noindent
\textbf{Proof} 
(Davis~\cite[8.7.1]{Dav})
\noindent
\textbf{QED}
\indent

From the Lemma we have
\[\mathrm{G}(x_1,x_2,x_3,\ldots,x_n) = A^{-1}\mathrm{G}(y_1,y_2,y_3,\ldots,y_n)\tilde{A^{-1}}\]
Then by the rules of the inversion of matrix products we have
\[\mathrm{G^{-1}}(x_1,x_2,x_3,\ldots,x_n) = \tilde{A}\mathrm{G^{-1}}(y_1,y_2,y_3,\ldots,y_n)A \]
Now the inner product of orthogonal polynomials gives a diagonal Gram matrix with
\[\lambda_i = (x_i,x_i)\]
So the inverse is 
\[\mathrm{G^{-1}}(y_1,y_2,y_3,\ldots,y_n) =  \left[
\begin{array}{cccc}
\lambda_1^{-1} & 0 & 0 & \ldots,0 \\ 
0 & \lambda_2^{-1} & 0 & \ldots,0 \\
0 & 0 & \lambda_3^{-1} & \ldots.0 \\
0 & (0 & 0 & \ldots,\lambda_n^{-1}
\end{array}
\right]
\]
Now the matrix A = $(a_{ij})$ is lower triangular while its conjugate transpose $\tilde{A}$ is upper triangular. This means that the inverse Gram matrix is given by 
\[\mathrm{G^{-1}}(x_1,x_2,x_3,\ldots,x_n) = \mathrm{(b_{ij})} = \sum_{k=j}^n  a_{ki}a_{kj} /\lambda_k \]
\section{The Reproducing Kernels}
\subsection{A Non-symmetric Orthogonal Polynomial}
\newcommand{\bin}[2]{
  \left(
    \begin{array}{@{}c@{}}
     #1 \\ #2
    \end{array}
   \right)
    }
The reproducing kernel K over a non-symmetric domain such as $(0,\infty)$ for Laguerre polynomials is given by 
\[\mathrm{K_{(n-1)}}(x,y) = \sum_{i=1}^n \sum_{j=1}^n b_{ij}x^{(i-1)}{y^{(j-1)}}\]
where
\[a_{ij} = \bin{i-1}{j-1} \frac{(-1)^{j-1}}{(j-1)!}\]
and
\[\lambda_i = 1\]
so
\[b_{ij} = \frac{(-1)^{i-1}}{(i-1)!} \frac{(-1)^{j-1}}{(j-1)!} \sum_{k=j}^n \bin{k-1}{i-1}\bin{k-1}{j-1}\]
The reproducing kernel K over symmetric domains such as $(-1,1)$ for Legendre polynomials and $(-\infty,-\infty)$ for Hermite polynomials naturally divides into separate even and odd kernels. The traditional polynomial expansion for these orthogonal polynomials is in descending order so that the even and odd terms can be included in one expression. However we require ascending polynomials which are more natural for our purposes.
\subsection{Even Orthogonal Polynomials}
The even reproducing kernel K is
\[\mathrm{K_{2(n-1)}}(x,y) = \sum_{i=1}^n \sum_{j=1}^n b_{ij}x^{2(i-1)}{y^{2(j-1)}}\]
For the Legendre case we have
\[a_{ij} = \frac{(-1)^{i+j}}{(i-1)!} \bin{i-1}{j-1} \frac{\Gamma(j-1/2+(i-1))}{\Gamma(j-1/2)}\]
and
\[\lambda_i = 2i-3/2\]
Therefore the inverse is given by
\[ b_{ij}  = \sum_{k=j}^n \frac{(-1)^{i+j}\bin{k-1}{i-1}\bin{k-1}{j-1}}{((k-1)!)^2(2k-3/2)} \frac{\Gamma(i-1/2+(k-1))}{\Gamma(i-1/2)} \frac{\Gamma(j-1/2+(k-1))}{\Gamma(j-1/2)} \]
For the Hermite case we have
\[a_{ij} = (2(i-1))! \frac{(-1)^{i-j}}{(i-j)!} \frac{2^{2(j-1)}}{(2(j-1))!}\]
and
\[\lambda_i = \Gamma(1/2) 2^{2(i-1)} (2(i-1))! \]
Therefore the inverse is given by
\[ b_{ij}  = \frac{(-1)^{i+j}}{\Gamma(1/2)} \frac{2^{2(i-1)}}{(2(i-1))!} \frac{2^{2(j-1)}}{(2(j-1))!} \sum_{k=j}^n \frac{(2(k-1))!}{2^{2(k-1)}} \frac{1}{(k-i)! (k-j)!} \]
\subsection{Odd Orthogonal Polynomials}
The odd reproducing kernel K is
\[\mathrm{K_{2(n-1)+1}}(x,y) = \sum_{i=1}^n \sum_{j=1}^n b_{ij}x^{2(i-1)+1}{y^{2(j-1)+1}}\]
For the Legendre case we have
\[a_{ij} = \frac{(-1)^{i+j}}{(i-1)!} \bin{i-1}{j-1} \frac{\Gamma(j+1/2+(i-1))}{\Gamma(j+1/2)}\]
and
\[\lambda_i = 2i-1/2\]
Therefore the inverse is given by
\[ b_{ij}  = \sum_{k=j}^n \frac{(-1)^{i+j}\bin{k-1}{i-1}\bin{k-1}{j-1}}{((k-1)!)^2( 2k-1/2)} \frac{\Gamma(i+1/2+(k-1))}{\Gamma(i+1/2)} \frac{\Gamma(j+1/2+(k-1))}{\Gamma(j+1/2)} \]
For the Hermite case we have
\[a_{ij} = (2(i-1)+1)! \frac{(-1)^{i-j}}{(i-j)!} \frac{2^{2(j-1)+1}}{(2(j-1)+1)!}\]
and
\[\lambda_i = \Gamma(1/2) 2^{2(i-1)+1} (2(i-1)+1)! \]
Therefore the inverse is given by
\[ b_{ij}  = \frac{(-1)^{i+j}}{\Gamma(1/2)} \frac{2^{2(i-1)+1}}{(2(i-1)+1)!} \frac{2^{2(j-1)+1}}{(2(j-1)+1)!} \sum_{k=j}^n \frac{(2(k-1)+1)!}{2^{2(k-1)+1}} \frac{1}{(k-i)! (k-j)!} \]
\subsection{Notes}
These are another equivalent form of the "kernel polynomials" referred to by $Szeg\ddot o$'s classic text on orthogonal polynomials.
That indicates that these are reproducing kernels for polynomials up to the same order as the size of the Gram matrix.
Also it is apparent that expectation of the kernel expansion is zero in the odd case and a tautology in the even or non-symmetric case.
The estimate of any function capable of integration gives the best possible polynomial estimate to the order of Gram matrix.\
\[f(x) \simeq \int_{\mathbf{D}}^{\ } K(x,y) f(y) dy \]
where $\mathbf{D}$ is the appropriate domain and K(x,y) is the appropriate kernel.
\section{Examples of condition numbers}
\newcommand{\pnorm}[2]{\|\, #1 \,\|_{#2}}
$\pnorm{A}{\infty}$ is the infinite $p$-norm which can be used to calculate the condition number 
\[\mathrm{\kappa (A_\infty) = \pnorm{A}{\infty}\pnorm{A}{\infty}^{-1}}\]
where
\[\mathrm{\pnorm{A}{\infty} = \max_{1 \le i\le n} \sum_{j=1}^n | a_{ij} | }\]
The details of the various classical orthogonal polynomials can be found in $Szeg\ddot o$.
\subsection{The Laguerre polynomial case}
\begin{center}
\begin{tabular}{ll}
Size & Condition Number\\
\hline
1 & 1 \\
2 & 9 \\
3 & 288 \\
4 & 22620 \\
5 & 2811960 \\
6 & 505744470 \\
7 & 125307922380 \\
8 & 48125908977898.5
\end{tabular}
\end{center}
As expected they grow rather dramatically!
\subsection{The Odd Legendre polynomial case}
\begin{center}
\begin{tabular}{ll}
Size & Condition Number\\
\hline
1 & 1 \\
2 & 37+1/3 \\
3 & 1251.375 \\
4 & 48319.7 \\
5 & 1627592.3610491073 \\
6 & 50651670.61979168 \\
7 & 1819734500.2385054 \\
8 & 62490018821.32725
\end{tabular}
\end{center}
As expected they grow rather dramatically!
\subsection{The Even Legendre polynomial case}
\begin{center}
\begin{tabular}{ll}
Size & Condition Number\\
\hline
1 & 1 \\
2 & 20 \\
3 & 513.1875 \\
4 & 18150 \\
5 & 607905.84765625 \\
6 & 18718753.030133925 \\
7 & 596835533.9163207 \\
8 & 20660541993.47005
\end{tabular}
\end{center}
As expected they grow rather dramatically!
\subsection{The Odd Hermite polynomial case}
\begin{center}
\begin{tabular}{ll}
Size & Condition Number\\
\hline
1 & 1 \\
2 & 18.375 \\
3 & 635.9765625 \\
4 & 63026.455078125 \\
5 & 9451396.300048828 \\
6 & 1970700695.5026627 \\
7 & 544334817099.9446 \\
8 & 192139661596606.84
\end{tabular}
\end{center}
As expected they grow rather dramatically!
\subsection{The Even Hermite polynomial case}
\begin{center}
\begin{tabular}{ll}
Size & Condition Number\\
\hline
1 & 1 \\
2 & 4.5 \\
3 & 114.84375 \\
4 & 6838.330078125 \\
5 & 699962.9919433594 \\
6 & 114340147.64968875 \\
7 & 31145186114.19374 \\
8 & 10775334775103.02
\end{tabular}
\end{center}
As expected they grow rather dramatically!
\section{Examples of error variance numbers}
\subsection{Trigonometric Functions over a finite domain}
\subsubsection{f(x) = sin(pi x)}
\begin{center}
\begin{tabular}{lll}
Size & Taylor Variance & Our Estimate Variance\\
\hline
2 & 0.80166669 & 0.00878023 \\
3 & 0.03778397 & 0.00003698 \\
4 & 0.00060558 & 4.90168598e-8 \\
5 & 4.21177985e-6 & 2.67561342e-11 \\
6 & 1.47725389e-8 & 7.097434670e-15 \\
7 & 2.89640460e-11 & 1.024073362e-18 \\
8 & 3.42305581e-14 & 8.722916936e-23
\end{tabular}
\end{center}
Our Estimate variances are significantly smaller than the Taylor equivalent variances which suggests that our Estimates could be very accurate when used as a replacement.
\subsubsection{f(x) = cos(pi x)}
\begin{center}
\begin{tabular}{lll}
Size & Taylor Variance & Our Estimate Variance\\
\hline
2 & 0.20346805 & 0.07606160 \\
3 & 0.00537462 & 0.00067401 \\
4 & 0.00005545 & 1.52457206e-6 \\
5 & 2.69759727e-7 & 1.26430443e-9 \\
6 & 6.99788566e-10 & 4.731342604e-13 \\
7 & 1.05658774e-12 & 9.147610639e-17 \\
8 & 9.91539162e-16 & 1.005215759e-20
\end{tabular}
\end{center}
Our Estimate variances are still significantly smaller than the Taylor equivalent variances which suggests that our Estimates could still be very accurate when used as a replacement.
Note that the Taylor variances here are significantly smaller than those  for the sine function. However our Estimate variances are larger than those for the sine function.
\subsection{Exponential Functions over a semi-infinite domain}
\subsubsection{f(x) = exp(-x)}
\begin{center}
\begin{tabular}{lll}
Size & Taylor Variance & Our Estimate Variance\\
\hline
2 & 5/6  & 1/48  \\
3 & 31/12 & 1/192 \\
4 & 209/24 & 1/768  \\
5 & 1471/48 & 1/3072  \\
6 & 10625/96  & 1/12288 \\
7 & 78079/192 & 1/49152 \\
8 & 580865/384 & 1/196608
\end{tabular}
\end{center}
The significant difference between the Taylor variance and our Estimate variance is due to the fact that the Taylor series is very good locally around its point of expansion while our polynomial estimate is very good more globally. This can be seen visually below for the case when the size of the Gram matrix is 8. The relevant polynomials are given below followed by their plot over the initial range.
Of course the larger the size of Gram matrix the better the approximation going further out from the origin.
\begin{eqnarray}
myEstimate&=&-1/1290240x^7 + 1/20480x^6 - 37/30720x^5 + 31/2048x^4 \nonumber\\ 
   & & {} - 163/1536x^3  + 219/512x^2 - 247/256x + 255/256 \nonumber \\
myTaylor&=&-1/5040x^7 + 1/720x^6 - 1/120x^5 + 1/24x^4 \nonumber\\
   & & {} - 1/6x^3 + 1/2x^2 - x +1 \nonumber
\end{eqnarray}
\sageplot{plot(-1/1290240*x^7 + 1/20480*x^6 - 37/30720*x^5 + 31/2048*x^4 - 163/1536*x^3 + 219/512*x^2 - 247/256*x + 255/256,0,3.5,ymin = -1/2,legend_label='myEstimate')
+plot(-1/5040*x^7 + 1/720*x^6 - 1/120*x^5 + 1/24*x^4 - 1/6*x^3 + 1/2*x^2 -x+1,0,3.5,ymin=-1/2,legend_label='myTaylor',color='red'), figsize=[4,2]}

This could mean that our estimate is a better way of approximating the Partition Function of Statistics Mechanics as well as the solution of linear constant coefficient Ordinary Differential Equations.

\section{Conclusion}
Expectations of our estimates are approximations to periods or exponential periods. In some cases this suggests some interesting L-series. This will be investigated further in subsequent work.
Since our estimates show how to carry out successful interpolation in 1D we can extend this to higher dimensions where the domain is a tensor product of 1D domains using Kronecker/Tensor Products. However as interpolation in higher dimensions in general is not always possible due to the frequent lack of unisolvence (Davis~\cite[2.4]{Dav}). There are various strategies which can be employed here which we leave to subsequent work.

\end{document}